\documentclass[11pt]{article}
\usepackage{amssymb,geompsfi,amscd}

\def\l{{\germ l}}

\def\C{{\mathbb C}}

\def\E{{\mathbb E}}

\def\H{{\mathbb H}}

\def\P{{\mathbb P}}
\def\Q{{\mathbb Q}}
\def\R{{\mathbb R}}

\def\Z{{\mathbb Z}}

\def\sL{{\cal L}}

\font\l=cmr10 at 10pt
\font\ls=cmr7
\font\lss=cmr5
\font\lsy=cmsy10
\font\lsys=cmsy7
\font\lsyss=cmsy5
\font\lmi=cmmi10
\font\lmis=cmmi7
\font\lmiss=cmmi5
\font\lex=cmex10

 at 6truept
 at 6truept

\def\mapright#1{\smash{
  \mathop{\longrightarrow}\limits^{#1}}}
\def\mapdown#1{\Big\downarrow
  \rlap{$\vcenter{\hbox{$\scriptstyle#1$}}$}}

\newtheorem{theorem}{Theorem}[section]
\newtheorem{lemma}[theorem]{Lemma}
\newtheorem{proposition}[theorem]{Proposition}

\newtheorem{corollary}[theorem]{Corollary}

\newcommand\heading[1]{\smallskip\noindent{\bf
#1}}
\newcommand\qed{\nopagebreak[4]\begin{flushright}\rule{0.1in}{0.1in}
\end{flushright}\pagebreak[2]}

\newcommand{\makefig}[3]{
	\begin{figure}[htbp]
        \refstepcounter{figure}
	\label{#2}
        \begin{center}
		~#3~\\
		\medskip
                {\sf Figure \thefigure.  #1}
        \end{center}
	\medskip
	\end{figure}
}

\def\GL{\mathrm{GL}}
\def\HH{\mathrm{H}}
\def\bdry{\partial}
\def\Star{\mathrm{Star}}

\title{Lehmer's Problem, McKay's Correspondence, and $2,3,7$
\smallskip
{\small Dedicated to the memory of Ruth Michler}}
\author{Eriko Hironaka}
\begin{document}
\maketitle

\section{Introduction}

This paper addresses a long standing open problem due to Lehmer
in which the triple 2,3,7 plays a notable role.  Lehmer's problem
asks whether there is a gap between 1 and the next smallest
algebraic integer with respect to Mahler measure.  The question
has been studied in a wide range of contexts including number theory,
ergodic theory, hyperbolic geometry, and knot theory;
and relates to basic questions such as describing
the distribution of heights of algebraic integers,
and of lengths of geodesics on arithmetic surfaces.
See, for example, \cite{Everest-Ward99} and \cite{G-H:Survey} 
for surveys and references.  This paper focuses on 
the role of Coxeter systems in Lehmer's problem.   The 
analysis also leads to a topological version of McKay's correspondence. 

We review some properties of Coxeter systems in Section 1, and 
Coxeter links in Section 2.
Section 3 covers Lehmer's problem, and Section 4 
contains some remarks on a topological generalization of 
McKay's correspondence.

\section{Coxeter Systems}

A {\it Coxeter system} consists of a vector space $V$ with a 
distinguished ordered basis $B = \{e_1,\dots,e_n\}$, and an inner product 
$$
\langle v_i,v_j \rangle = -2\cos{\frac{\pi}{m_{i,j}}}, 
$$
where $m_{i,i} = 1$, and if $i\neq j$, 
$m_{i,j} \in \{2,3,\dots,\infty\}$.  Associated to the Coxeter system
is the {\it Coxeter group} $G \subset \GL(V)$ generated by
reflections $S=\{s_1,\dots,s_n\}$ through hyperplanes perpendicular to 
$e_1,\dots,e_n$ respectively.
The action of $s_i \in S$ is given by
$$
s_i(e_j) = e_j - \langle e_i,e_j \rangle e_i.
$$
The group  $G$ has presentation  
$$
G = \langle\ s_1,\dots,s_n\ : \ (s_is_j)^{m_{i,j}} = 1\ \rangle.
$$
Coxeter systems are typically denoted by $(G,S)$.

A Coxeter system is determined by its {\it Coxeter graph} $\Gamma$.
This is the graph 
with vertices $\nu_1,\dots,\nu_n$ corresponding to the elements of 
$S$ and edges labeled $m_{i,j}$ connecting distinct vertices $\nu_i$ and 
$\nu_j$ whenever $m_{i,j} > 2$. 

The {\it Coxeter element} of $(G,S)$ is the product of reflections
$$C = s_1\dots s_n,$$  
and is an important invariant of the system.  For
example, the spectral radius of the Coxeter element equals 1
if and only if $G$ is spherical or affine \cite{ACampo:Coxeter},
\cite{Howlett:Coxeter}.    

\subsection{Coxeter Links}\label{link-section}

A link $K$ in $S^3$ is {\it fibered} with fiber $\Sigma$ if
$\bdry \Sigma = K$ and $S^3 \setminus \Sigma = \Sigma \times I$.
If $K$ is a fibered link, there is an associated fibration
$$
\begin{array}{ccc}
\Sigma &\mapright{} &S^3 \setminus K \\
&&\mapdown{} \\
&&S^1
\end{array}
$$
The gluing map $h: \Sigma \rightarrow \Sigma$ induces a
{\it monodromy}:  
$$h_* : \mbox{H}_1(\Sigma; \R) \rightarrow \mbox{H}_1(\Sigma; \R).
$$
The {\it Alexander polynomial} of $K$ is the characteristic polynomial 
of $h_*$.

In this section we will only deal with {\it simply-laced} Coxeter systems,
where $m_{i,j} \in \{1,2,3\}$.  Since in this case all edges on
the Coxeter graph are labeled 3, we drop the labeling. 

\makefig{}{cross}
{\psfig{figure=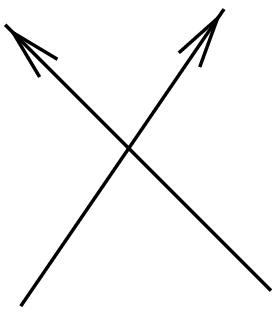,height=0.7in}}

An {\it ordered chord system} is a collection of oriented chords
embedded on a disk $D$ so that the endpoints lie on the boundary of
$D$.   A chord system 
$$
\{\ell_1,\dots,\ell_n\}; 
$$
is {\it positive} if $\ell_i$ intersects $\ell_j$ positively 
whenever $i>j$.
Figure~\ref{cross} shows two chords segments intersecting
positively.  

Let $A$ be the intersection matrix of the chords in $\sL$ (with diagonal
entries equal to zero.)  Then $\sL$ is a positive chord system if and only
if the lower triangular part of $A$ is non-negative (and hence the upper
triangular part is non-positive.)

From an ordered chord system $\sL$, we can define a fibered link $K_{\sL}$ 
with fiber $\Sigma_{\sL}$ as follows.
Consider $\sL$ as an ordered collection of chords
on a disk $D$ embedded as the unit disk in the $x,y$ plane in $\R^3$.  
Attach bands to $D$ with one full positive twist as in 
Figure~\ref{mur} in the order 
given by the ordering on the chord system.  That is, the twisted
band $\eta_i$ corresponding to the chord $\ell_i$ lies over 
the band
$\eta_j$ corresponding to the chord $\ell_j$ if and only if $i > j$.
This defines a surface $\Sigma_{\sL} \subset \R^3$ and its link
boundary $K_{\sL} = \bdry\Sigma_{\sL}$.
Identifying $\R^3$ with the subset of $S^3$ considered as the one
point compactification of $\R^3$ gives a link $K_{\sL}$ with 
specific Seifert surface $\Sigma_{\sL}$.

\makefig{Murasugi sum.}{mur}
{\psfig{figure=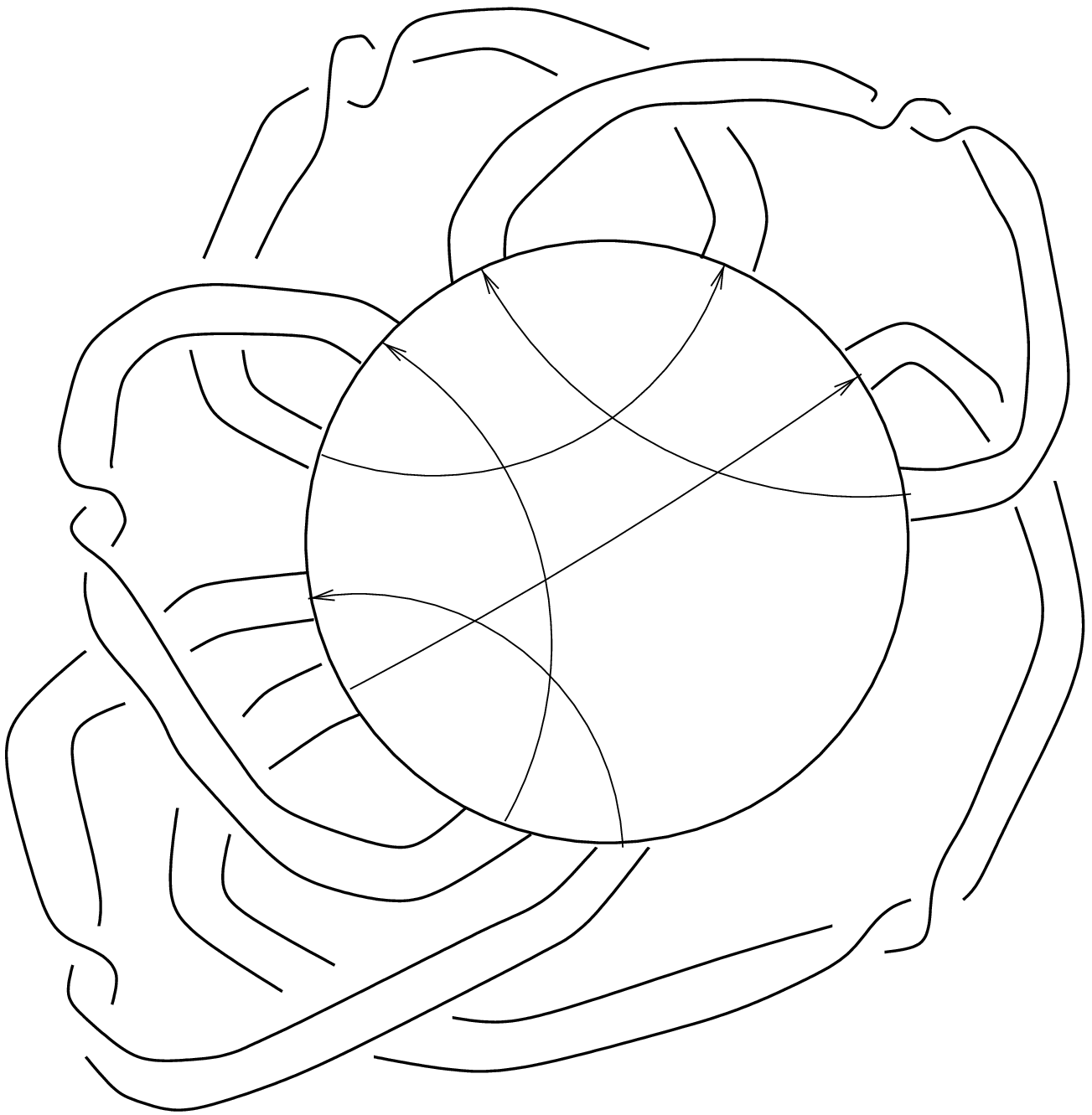,height=1.35in}}

Stallings showed (see \cite{Gabai:Murasugi}) that the
Murasugi sum of two fibered links is fibered.  Thus, the links
$K_{\sL}$ are fibered with fiber $\Sigma_{\sL}$.  The oriented chords
in $\sL$ extend to give a basis for $\HH_1(\Sigma_{\sL};\R)$.
Let $\Gamma$ be the incidence graph of $\sL$.
If $\sL$ is a positive chord system, then we call the pair
$(K_{\sL},\Sigma_{\sL})$ a {\it Coxeter link} associated to $(G,S)$.
Its Seifert matrix equals $M = I + A^+$ and hence $h_* = M^tM^{-1}$,
where $A$ is the intersection matrix of the chord system
defined above.  The bilinear form of
$(G,S)$ can be written as
$$
B = M + M^t
$$
and the Coxeter element of $(G,S)$ equals 
(cf \cite{Howlett:Coxeter})
$$
C = -M^tM^{-1} = -h_*.
$$

\subsection{Realizable graphs}\label{realization-section}

A graph $\Gamma$ is {\it realizable} if it is the incidence
graph of a chord diagram.  Here are some examples of realizable graphs:

\begin{description}
\item{(i)} Complete graphs;
\item{(ii)} Cyclic graphs;
\item{(iii)} Join of two realizable graphs at 
one vertex; and
\item{(iv)} Trees.
\end{description}

There are, however, obstructions to realizability.  Figure~\ref{cube}
gives an example of a non-realizable graph.
Let $\Gamma$ be a graph with vertices $S$.  A subgraph 
$\Gamma' \subset \Gamma$ is an {\it induced subgraph} if for
some $S' \subset S$, $\Gamma'$ is the subgraph containing all 
edges on $\Gamma$ whose endpoints are in $S'$. An induced cycle in 
$\Gamma$ is a cycle which is an induced subgraph.

\begin{proposition}$\Gamma$ is not realizable if there is a subset $S' \subset S$
such that
\begin{description}
\item {(i)} $S'$ contains at least three vertices;
\item {(ii)} $S'$ is disjoint;
\item {(iii)} there is an $s \in S$ so that $s$ is joined by an edge in $\Gamma$
to every vertex in $S'$; and
\item {(iv)} there is an induced cycle in $\Gamma$ containing $S'$.
\end{description}
\end{proposition}

\makefig{Non-realizable graph.}{cube}
{\psfig{figure=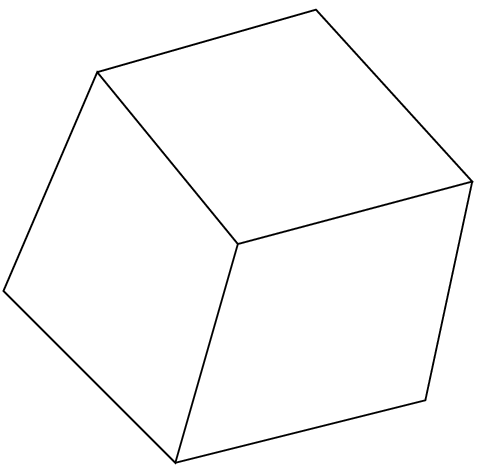,height=0.5in}}

As was pointed out to me by R. Vogeler, any graph can be realized as an 
incidence graph if we
generalize to higher dimensional diagrams.  A higher dimensional chord
system is a union of mutually transverse, linear embeddings 
$$
(D,\bdry D) \hookrightarrow (B,\bdry B)
$$ 
where $D$ is an $n$-disk and $B$ is an $n+1$-ball embedded in $\R^{n+2}$.
\begin{lemma}  Any finite graph with no self-loops or double edges
can be realized as the incidence graph of some higher dimensional chord
diagram.
\end{lemma}
An analysis of higher dimensional Coxeter links will be the topic of a
future article.

\subsection{Comments on ordering and positivity}\label{ordering-section}

In Section~\ref{link-section}, we saw that a Coxeter graph $\Gamma$
has a corresponding Coxeter link, if it is both realizable as a 
chord system, and its ordering is compatible with a positive ordering
on the chord system.  

\begin{lemma} Any chord diagram admits an ordering and orientation 
which is positive.
\end{lemma}

\heading{Proof.} Choose a direction vector $v$ from the 
center of the disk and orient the chords so that their 
direction vectors have positive inner product in the usual Euclidean
metric on $\R^2$ with $v$.  Now order the chords counter-clockwise
starting with the chord pointing furthest to the right of $v$.\qed

Not all orderings on a realizable graph, however, are realizable by
a positive chord system, see for example
Figure~\ref{squareorder}.
\makefig{Ordered graph which cannot be realized by a positive chord system}{squareorder}
{\psfig{figure=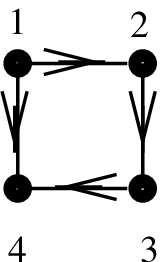,height=0.75in}}

As the examples in Section~\ref{examples-section} show, different
orderings on a chord system can give rise to different links.
To determine the link, however, it is not necessary to have all
the information of the ordering.

Given an ordered graph, there is an associated directed graph,
where edges are directed so that they point to the vertex with
larger index.  As one can see from the construction, we 
do not need all the information of the ordering on the chord diagram.
\begin{lemma}\label{ordering-lemma} If two orderings on a chord 
diagram $\sL$ have the same directed incidence graph then the
resulting fibered links are the same.
\end{lemma}
Lemma~\ref{ordering-lemma} is analogous to a result in Coxeter graph 
theory, which states that Coxeter elements depend only on
the directed Coxeter graph \cite{Shi:Enumeration}.  
\begin{lemma} The Coxeter element of Coxeter system $(G,S)$
depends only on the directed Coxeter graph $\Gamma$ of $(G,S)$.
\end{lemma}

The case when the incidence graph of $\sL$ is a tree has been
well studied, and the corresponding link has been called an
{\it arborescent link} \cite{Conway:Knots}.  Arborescent links
also appear as  {\it slalom links} in \cite{ACampo:Slalom}. 
Since any tree is realizable, there exists a Coxeter link associated
to any tree.  The corresponding Coxeter link will be independent of
ordering, but will depend on the realization.  Thus, for example, it is
possible to find non-equivalent links which are Coxeter
links for the same (ordered) Coxeter system.
\begin{lemma} If $\Gamma$ is an (unordered) tree with nodes of degrees 
less than or equal to 3, then the Coxeter link associated to $\Gamma$ 
is uniquely determined by $\Gamma$.
\end{lemma}
\makefig{Coxeter Link for a star graph.}{pretzel}
{\psfig{figure=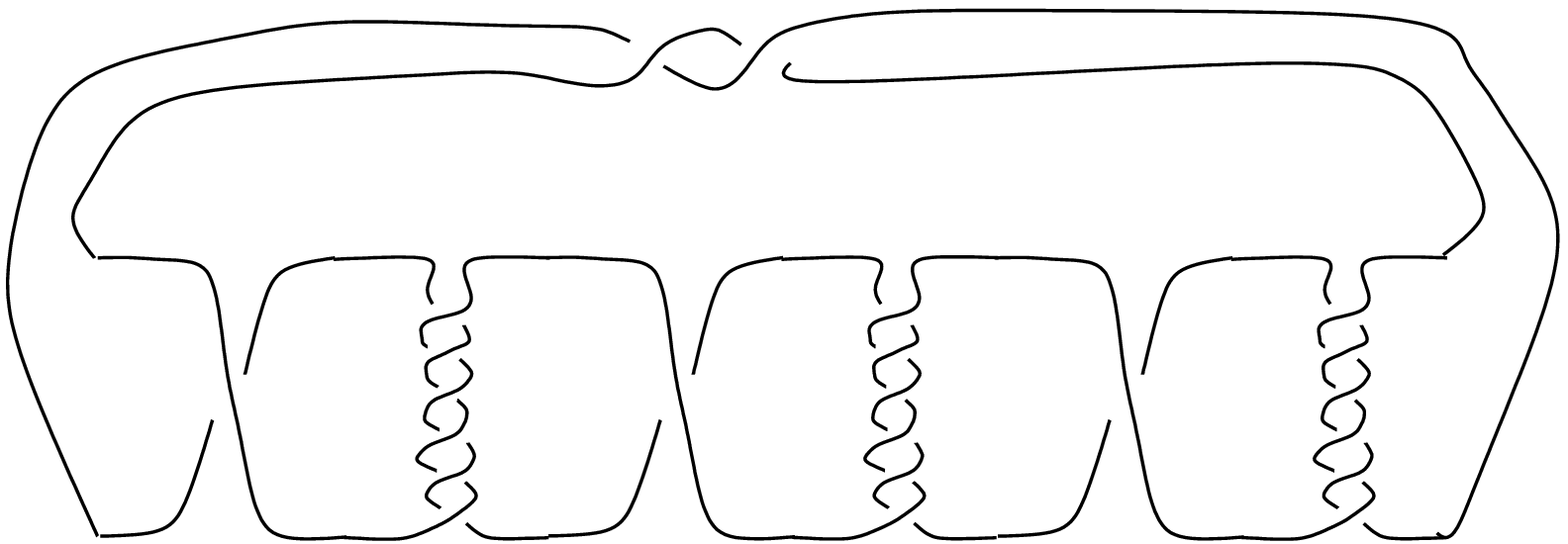,height=0.75in}}

\subsection{Examples}\label{examples-section}

\makefig{Star diagram and embedding.}{star}
{\psfig{figure=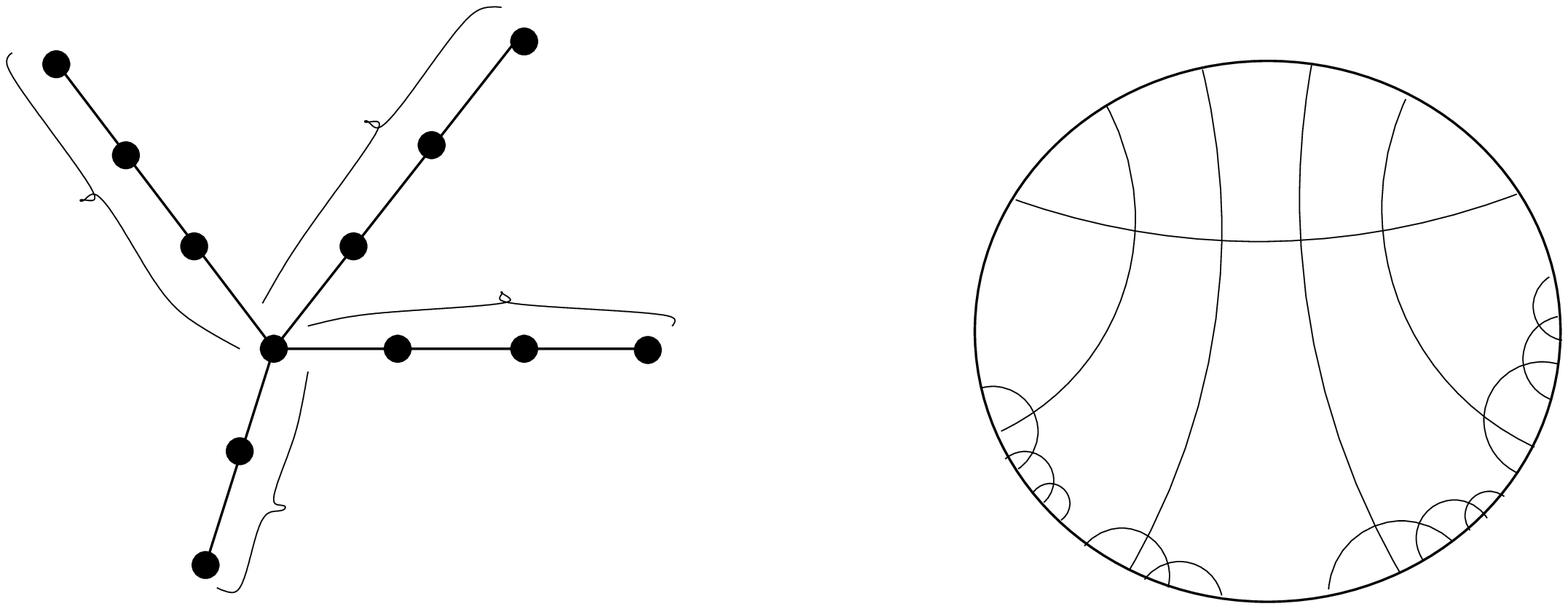,height=1in}}
We begin with the classical Dynkin diagrams.  
using a correspondence between {\it star-diagrams}
shown in Figure~\ref{star} 
and the $K_{p_1,\dots,p_k}$ pretzel links 
shown in Figure~\ref{pretzel}.
Since Dynkin diagrams are
trees with vertices of degree at most three, the associated link 
doesn't depend on their realization as a chord system.  
The correspondence Dynkin diagrams and their corresponding
Coxeter links is shown in Figure~\ref{Klein}.
\makefig{Links associated to Dynkin diagrams.}{Klein}
{\psfig{figure=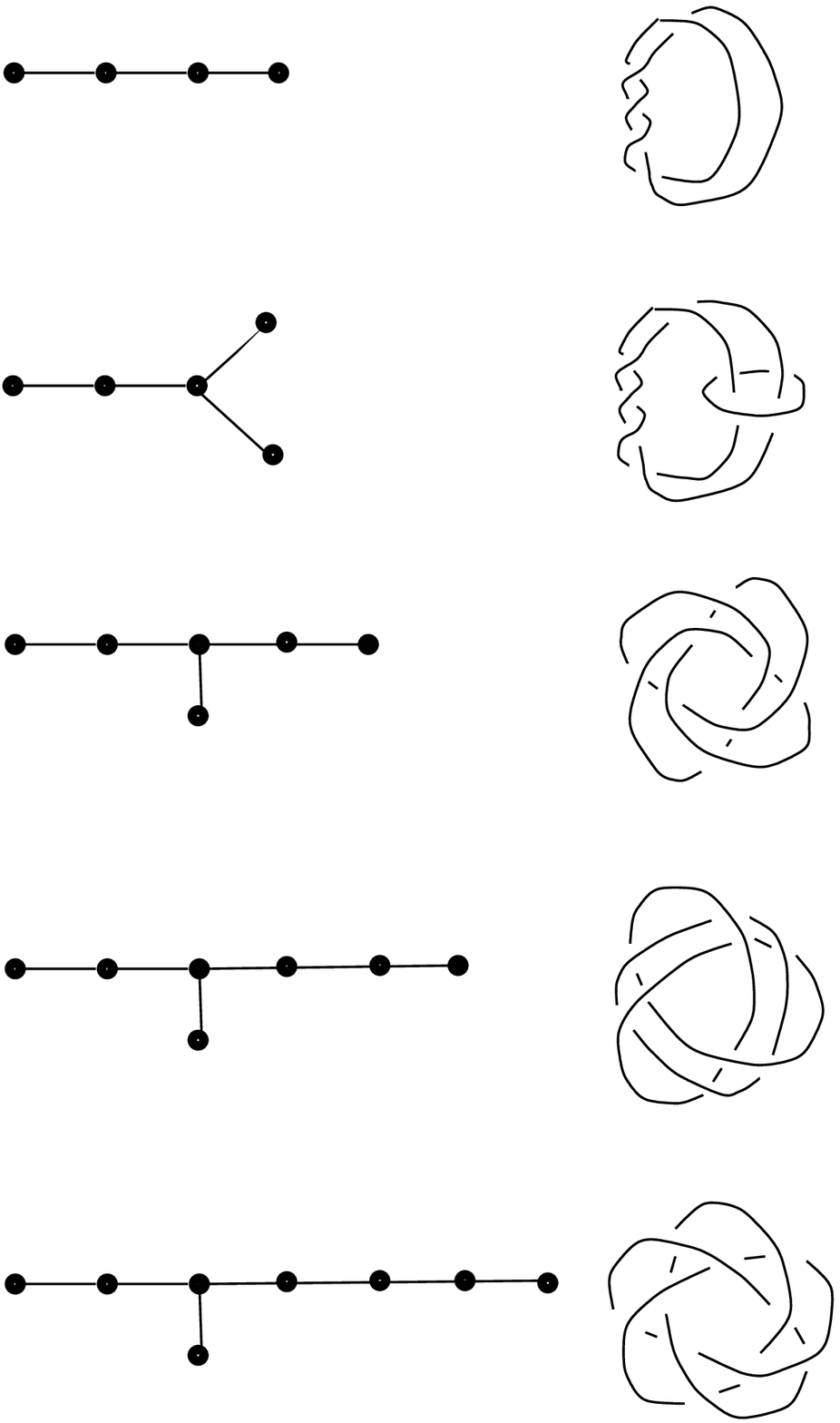,width=2in}}

It is possible to cook up examples of non-equivalent
links associated to the same ordered Coxeter system using
star-diagrams.  Take the two realization of the same tree 
shown in Figure~\ref{treeorder}.
One sees that the link on the left has two knotted components,
while the one on the right has a component which is the unknot,
hence the knots are not isotopic equivalent.
\makefig{Two embeddings of the same tree.}{treeorder}
{\psfig{figure=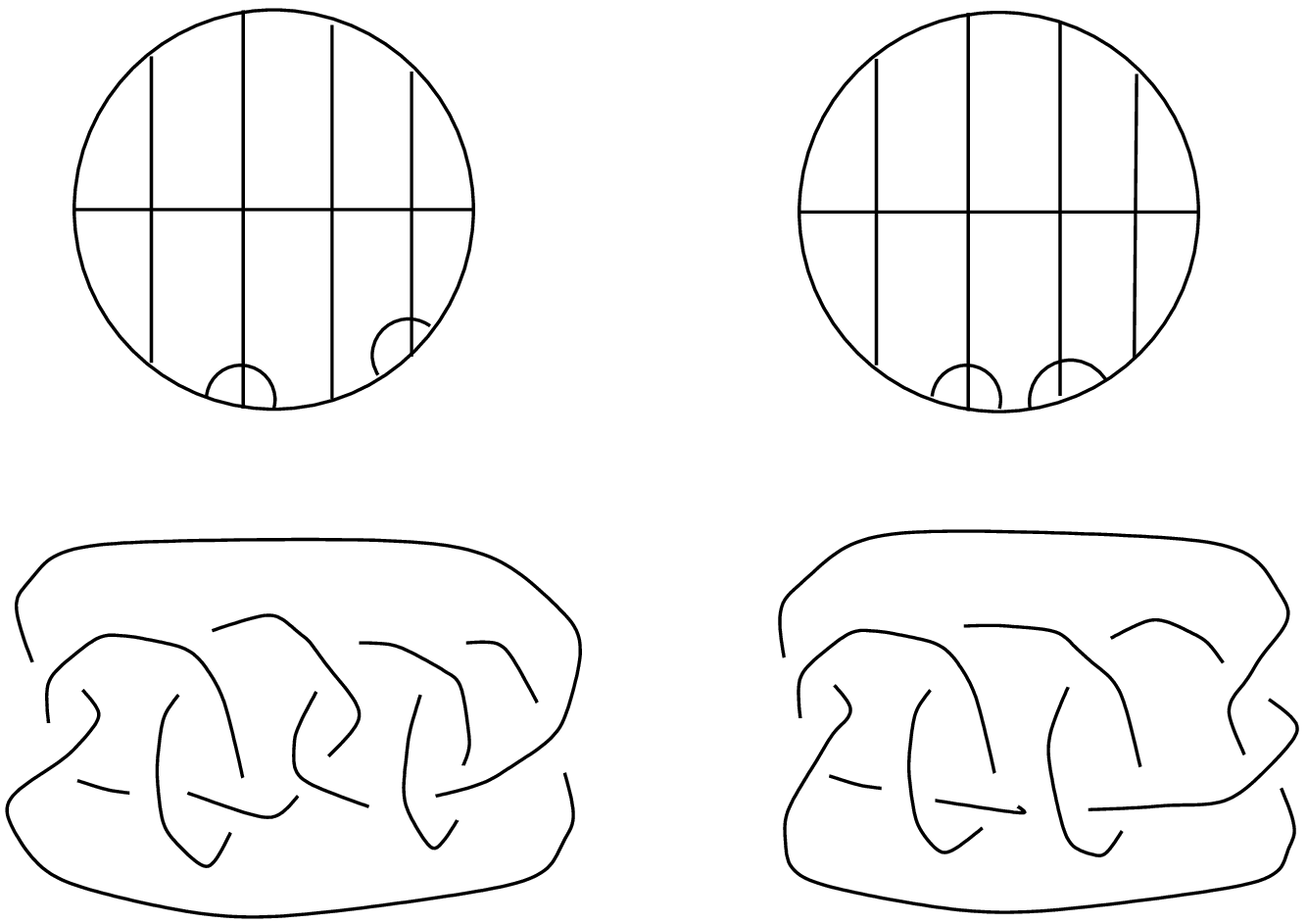,height=1.75in}}

When the incidence graph of a chord system contains cycles 
different orderings on the graph can give rise to different links. 
Consider for example, the 5-cycle.  Up to isotopy, there is only 
one chord diagram with this incidence graph, but there are 
two inequivalent positive orderings as shown in Figure~\ref{5-cycle-graph}.
\makefig{Two orientations for the 5-cycle}{5-cycle-graph}
{\psfig{figure=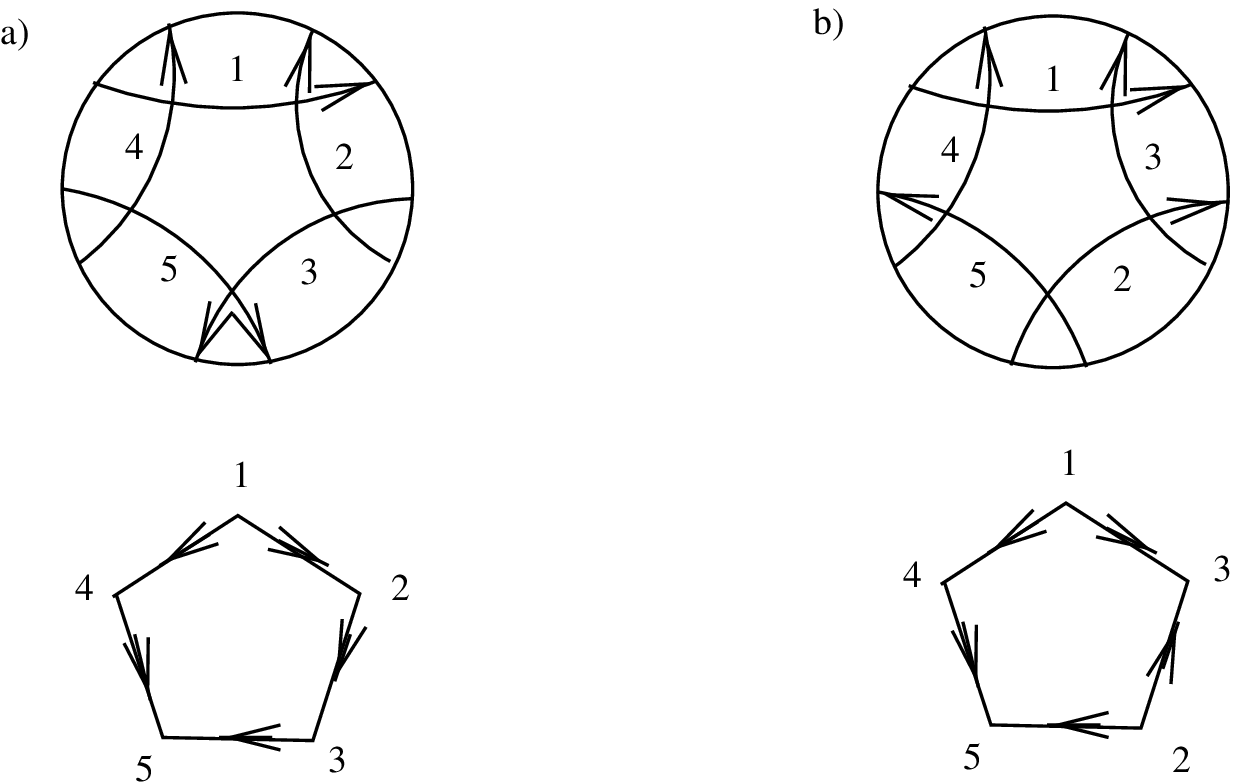,height=2in}}
The different orientations give rise to distinct characteristic
polynomials for the Coxeter elements:
\begin{eqnarray*}
\Delta_1(-t) &=& 1 - t - t^4 + t^5;\quad\mbox{and}\\
\Delta_2(-t) &=& 1 - t^2 - t^3 + t^5.\\
\end{eqnarray*}
These are computed using the Seifert matrices
$$
\left (
\begin{array}{ccccc}
1 & -1 & 0 & 0 & -1\\
0 & 1 & -1 & 0 & 0 \\
0 & 0 & 1 & -1 & 0\\
0 & 0 & 0 & 1 & -1\\
0 & 0 & 0 & 0 & 1
\end{array}
\right )
\qquad
\left (
\begin{array}{ccccc}
1 & 0 & -1 & -1 & 0\\
0 & 1 & -1 & 0 & -1 \\
0 & 0 & 1 & 0 & 0\\
0 & 0 & 0 & 1 & -1\\
0 & 0 & 0 & 0 & 1
\end{array}
\right ).
$$
The corresponding links are given in Figure~\ref{5-cycle-links}.
\makefig{Two Coxeter links for the 5-cycle}{5-cycle-links}
{\psfig{figure=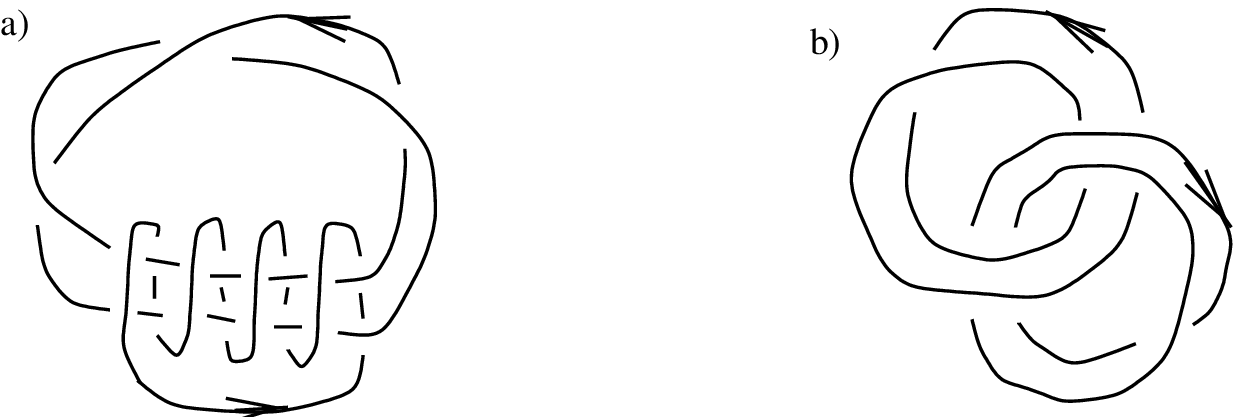,height=1in}}
They are distinct iterated torus links.

The simply-laced minimal hyperbolic Coxeter system
of smallest dimension is a triangle with a tail.
The Coxeter link (see Figure~\ref{hypert}) is uniquely determined
in this case by the requirement of positivity, and equals the 
mirror of the $10{}_{145}$-knot in Rolfsen's table \cite{Rolfsen76}
($22,3,3-$ in Conway's notation \cite{Conway:Knots}.) 
\makefig{Coxeter link associated to smallest hyperbolic Coxeter system}{hypert}
{\psfig{figure=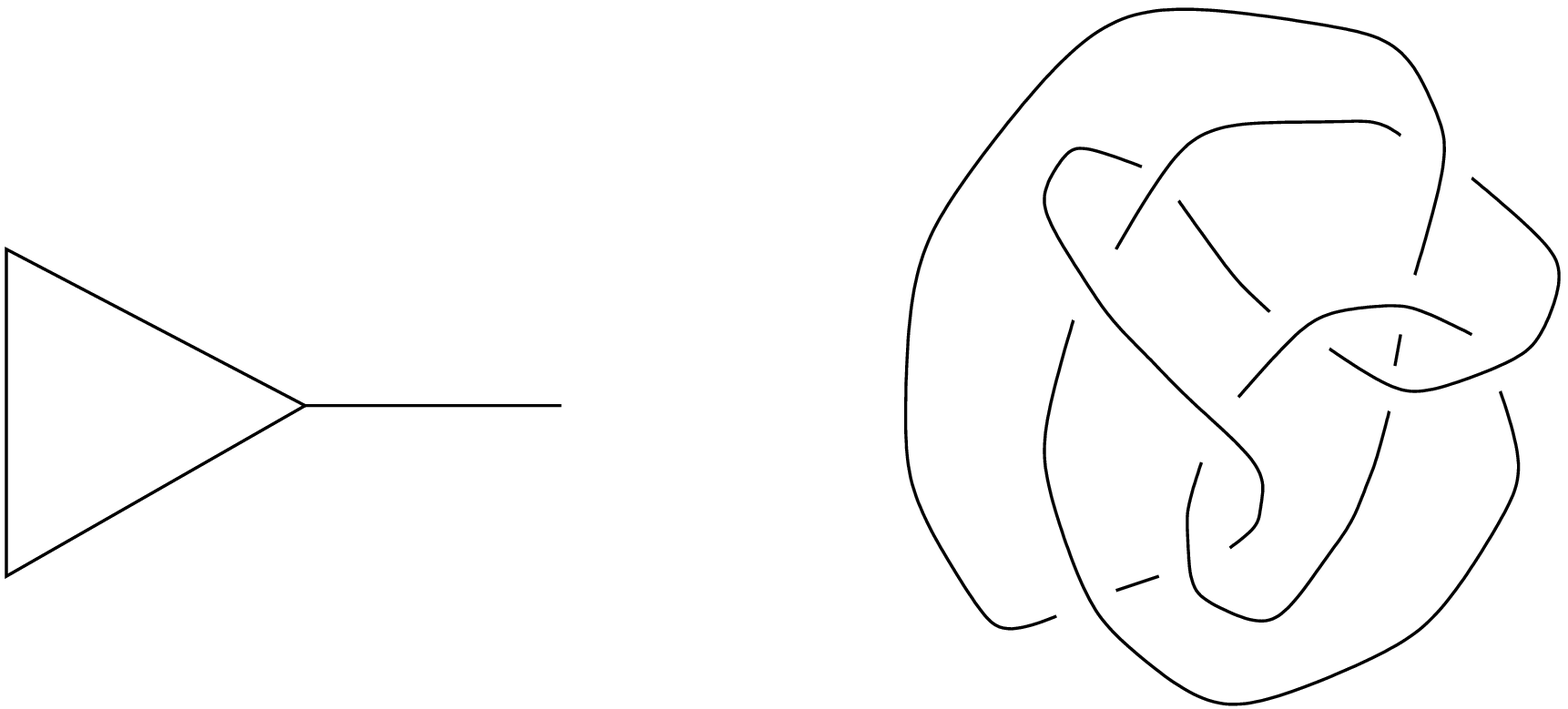,height=1in}}

\subsection{Remarks on the Geometry of Coxeter systems and Coxeter links}

A Coxeter system is {\it spherical} if its Coxeter group is a finite
Euclidean reflection group (reflections across hyperplanes through
the origin).  In the simply-laced case, these are 
the $A$-$D$-$E$ Coxeter systems.

A Coxeter system is {\it affine} if its Coxeter group
is an affine reflection group.  The simply laced ones are the
$\widetilde A$-$\widetilde D$-$\widetilde E$ Coxeter systems.  

The following is a classical fact due to Bourbaki \cite{Bourbaki:Lie}.

\begin{theorem} Let $\Gamma$ be a Coxeter system, and let
$B$ be its associated bilinear form.
Then $\Gamma$ is spherical if and only if $B$ is positive definite,
and affine if and only if $B$ is positive semi-definite.
\end{theorem}

A'Campo \cite{ACampo:Coxeter} and  Howlett \cite{Howlett:Coxeter} proved
the following relation between the geometry of the Coxeter system and
the Coxeter element.   

\begin{theorem}
Let $\Gamma$ be a Coxeter 
system and $C$ a Coxeter element.  Then
\begin{itemize}
\item The eigenvalues of $C$ lie on $\R \cup S^1$.
\item $\Gamma$ is spherical if and only if $C$ has finite order; 
\item $\Gamma$ is neither affine or spherical if and only if $C$ has 
an eigenvalue greater than one.
\end{itemize}
\end{theorem}

Links $K$ in $S^3$ have an analogous classification as torus links and 
satellite links, which include all algebraic links, 
and hyperbolic links, whose complements are uniformized by the
$3$-ball.  The monodromy $h_*$ has finite order and hence its
spectral radius is $1$ if $K$ is a fibered iterated torus link,
for example, when $K$ is an algebraic link.   

The above discussions bring up the following questions:
Let $K_\sL$ be the Coxeter link associated to a simply laced 
Coxeter system $(G,S)$.  Is it true that $K_\sL$ is hyperbolic 
if and only if the $(G,S)$ has indefinite intersection form?
When is the monodromy $h$ pseudo-Anosov?

\section{Lehmer's Problem}

Let $\alpha$ be an algebraic integer, and define its {\it size} to be
$$
\| \alpha \| = \prod_{\beta = \alpha^\sigma, |\beta|>1} \beta.
$$
This is also known as the {\it Mahler measure} of the minimal polynomial of
$\alpha$ over $\Q$.  It is well known that 
$\|\alpha\| = 1$ if and only if $\alpha^N = 1$ some $N$.
Thus we are interested in algebraic integers which are not roots of unity.

In 1933, Lehmer \cite{Lehmer33} asks 
whether for each $\delta > 1$, there exists
an algebraic integer such that 
$$
1 < \| \alpha \| < 1 + \delta.
$$
It is an easy exercise to see that the statement is false if one fixes
the degree of the integer.

Lehmer found polynomials with smallest Mahler measure for small degrees
and states in \cite{Lehmer33} that the smallest he could find of 
degree 10 has minimal polynomial
$$
P_L(x) = x^{10} + x^9 - x^7 - x^6 - x^5 - x^4 - x^3 + x + 1.
$$
Boyd \cite{Boyd89} and Mossinghoff \cite{Mossinghoff98} have done
searches up to degree 40, but so far no one has found a monic noncyclotomic
integer polynomial with smaller Mahler measure.

\makefig{Roots of the Lehmer polynomial}{Lehmer-roots}
{\psfig{figure=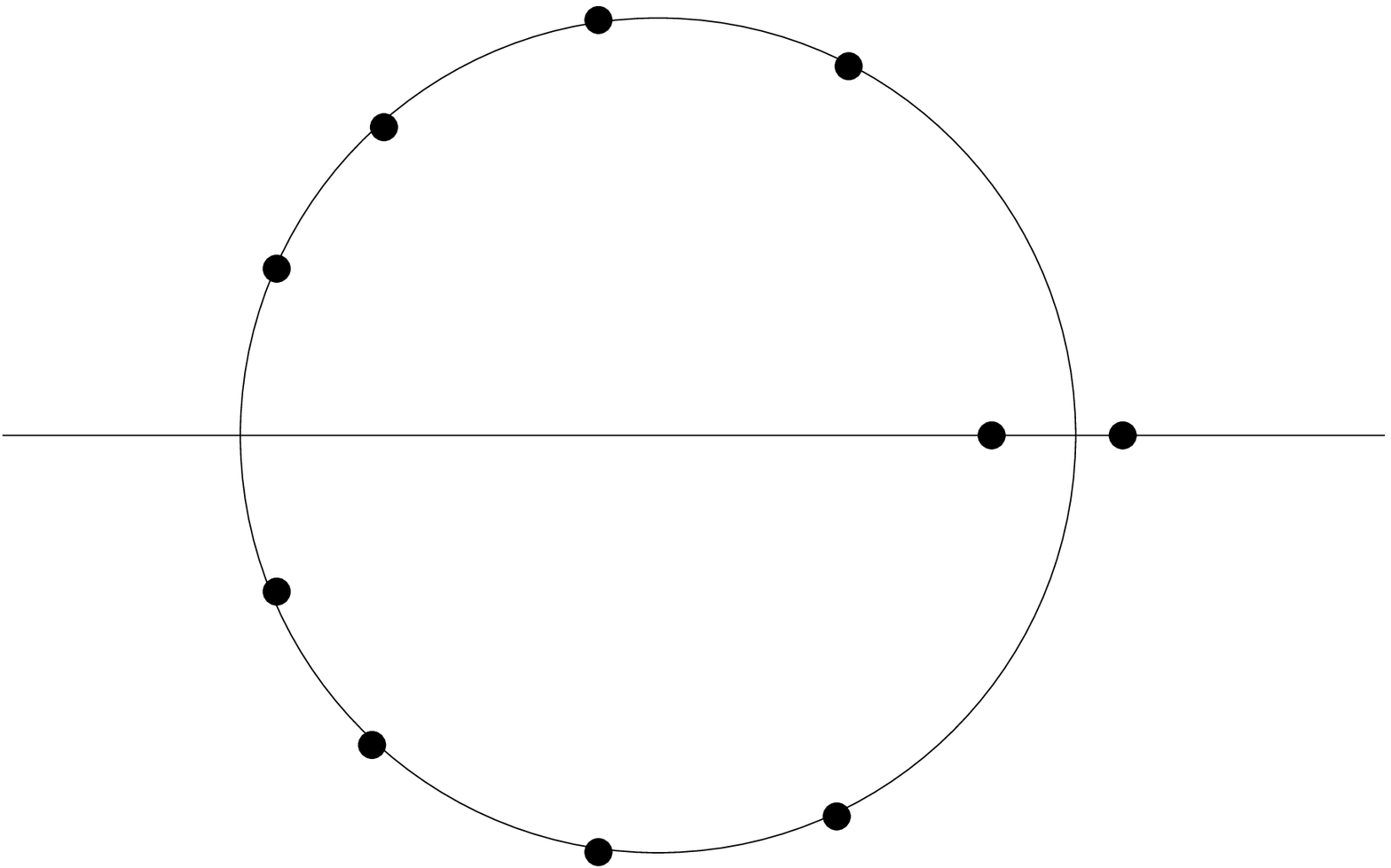,height=1in}}

One observes immediately that $P_L(x)$ is reciprocal, that is 
$$
P_L(x) = x^{d}P_L(\frac{1}{x}),
$$
where $d$ is the degree of $P_L(x)$ (in this case $d=10$.)
From Figure~\ref{Lehmer-roots}, one sees also that $P_L(x)$ has only one
root $\alpha_L = 1.17628\dots$ which we will call {\it Lehmer's number}
outside the unit circle.  Thus, $\alpha_L$ is what is known 
as a {\it Salem number}, that is, an algebraic integer
whose algebraic conjugates lie on or within the unit circle, with at 
least one conjugate on the unit circle (making the minimal polynomial
necessarily reciprocal.)   It is not known whether
there exist Salem numbers smaller than $\alpha_L$.

Smyth \cite{Smyth70} shows that among non-reciprocal polynomials
Lehmer's statement is false, and the
smallest Mahler measure $1.32472\dots$ is attained by
$$
x^3 - x + 1.
$$
Thus, it remains to determine whether there is a similar minimum for
Mahler measures of reciprocal monic integer polynomials.

It has been observed in various contexts that Lehmer's problem 
is related to the triple $(2,3,7)$ and more abstractly to the notion
of minimal hyperbolicity.  Before  going to examples, 
it is worth remarking that the triple has the simple distinguishing 
property that, among all $k$-tuples of positive integers $(p_1,\dots,p_k)$, 
$(2,3,7)$ gives the minimal positive value for 
$$
k-2 - \sum_{i=1}^k\frac{1}{p_k}. 
$$
This property comes into play in the minimality of Lehmer's
number $\alpha_L$ among the series of Salem numbers and algebraic numbers
which we describe in this section.

\subsection{Growth rates and the ($2,3,7$)-Triangle Group}\label{growth-section}

\makefig{Tiling of the hyperbolic disk by the action of $T_{2,3,7}$.}{t237}
{\psfig{figure=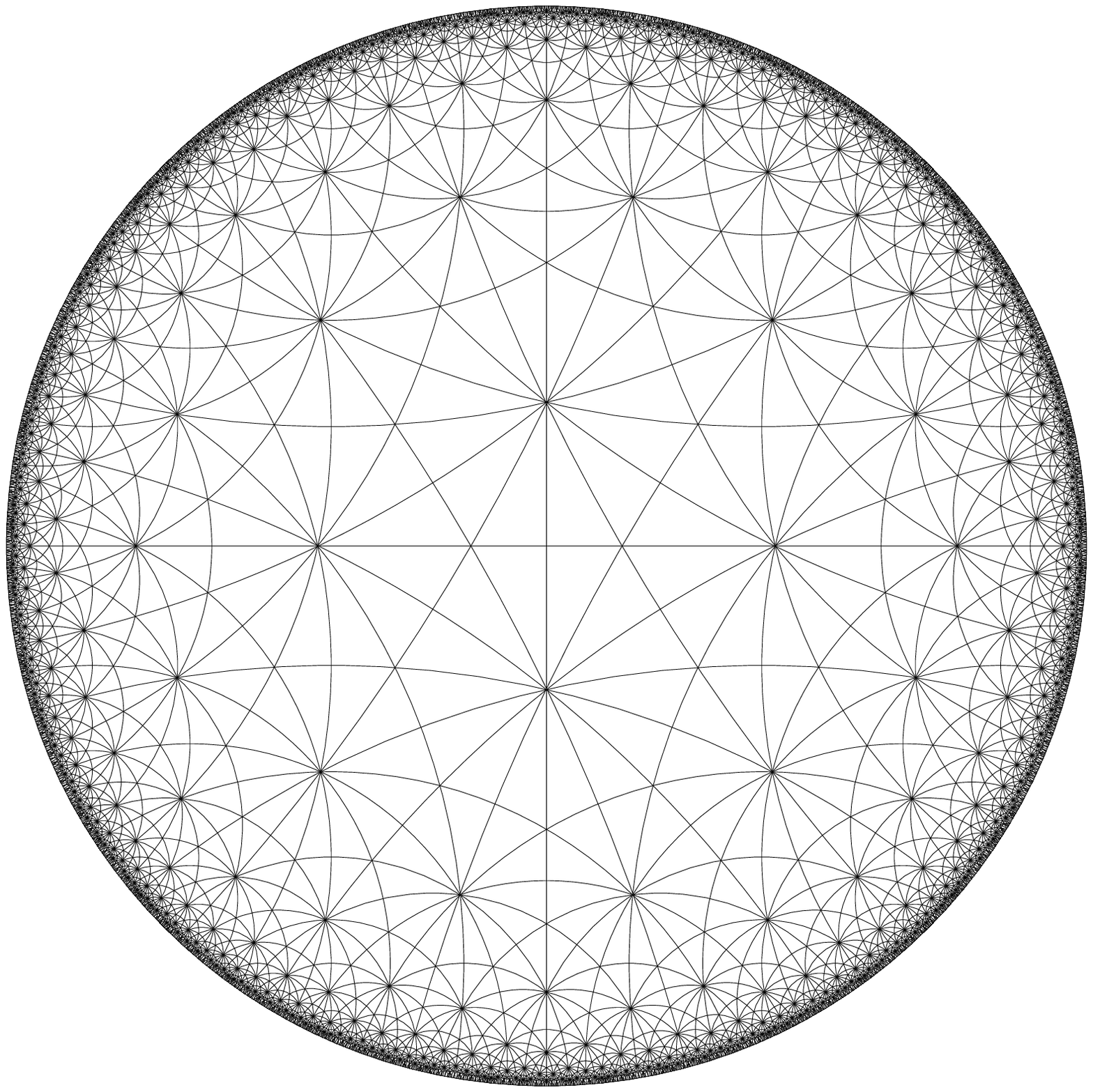,height=1.75in}}

Consider any pair $(G,S)$, where $G$ is a group and $S$ is a set of generators.
Let
$$
w_n = \mbox{number of words of minimal word length $n$ in $S$}
$$

The {\it growth series} of $(G,S)$ is the formal power series
$$
f_{(G,S)} = \sum_{n=1}^{\infty} w_n t^n
$$
and the growth rate $\alpha$ equals
$$\alpha = \frac{1}{\mbox{radius of convergence of $f_{(G,S)}$}}.$$
Another way to say this is that $w_n$ grows like $\alpha^n$ as $n$ 
gets large.

Let $G=T_{p_1,\dots,p_k}$ be a polygonal reflection group acting
on $S^2$, $\E^2$ or $\H^2$.   The group $G$ is the Coxeter group 
associated to the cyclic Coxeter graph with edges labeled $p_1,\dots,p_k$, 
and has presentation 
$$
G=\langle\ s_1,\dots,s_k\ : \ (s_i s_j)^{m_{i,j}}\ \rangle
$$
In this case Floyd  and Plotnick \cite{Floyd-Plotnick88}
show the following. (See also, \cite{Cannon-Wagreich92}.)  

\begin{theorem}  The growth series for $(G,S)$ is rational
$$
f_{(G,S)} = \frac{R(t)}{\Delta_{p_1,\dots,p_k}(t)}
$$
where $\Delta_{p_1,\dots,p_k}(t)$
is given by:
$$
\Delta_{p_1,\dots,p_k}(t) = (x - k + 1)\prod_{i=1}^k [p_i]\ +\ \sum_{i=1}^k
[p_1]\dots\widehat{[p_i]}\dots[p_k],  
$$
where $[p] = 1 + x + \dots + x^{p-1}$.
Furthermore, $\Delta_{p_1,\dots,p_k}(t)$ is a reciprocal monic integer
polynomial with a root, necessarily a Salem number, 
outside the unit circle if and only if 
$$
\frac{1}{p_1} + \cdots + \frac{1}{p_k} < k-2.
$$
\end{theorem}

\noindent
Thus, the growth rate $\alpha_{p_1,\dots,p_k}$  of $(G,S)$ is a 
Salem number if and only if 
$$
\chi  = \frac{1}{p_1} + \cdots + \frac{1}{p_k} - (k-2) < 0.  
$$
The polynomial $\Delta_{2,3,7}(x)$ equals Lehmer's polynomial $P_L(x)$, and
hence the triangle group $T_{2,3,7}$ has growth rate equal to $\alpha_L$. 

For the family of polynomials $\Delta_{p_1,\dots,p_k}(x)$ Lehmer's problem
is solved \cite{Hironaka98}. 

\begin{theorem}\label{delta-theorem}
Among the polynomials $\Delta_{p_1,\dots,p_k}(x)$, the one with 
smallest Mahler measure is Lehmer's polynomial $\Delta_{2,3,7}(x)$.
\end{theorem}

\noindent
This result is suggestive since among hyperbolic orbifold spheres
$(S^2;p_1,\dots,p_k)$, the one with maximal orbifold Euler characteristic
$\chi$ and minimal hyperbolic area is $(S^2;2,3,7)$.

\subsection {Alexander polynomials and the $(2,3,7,-1)$ Pretzel Knot}

Reidemeister \cite{Reidemeister32} remarked that the $(-2,3,7)$-pretzel
knot shown in Figure~\ref{p237} has Alexander polynomial $P_L(-x)$.

\makefig{}{p237}
{\psfig{figure=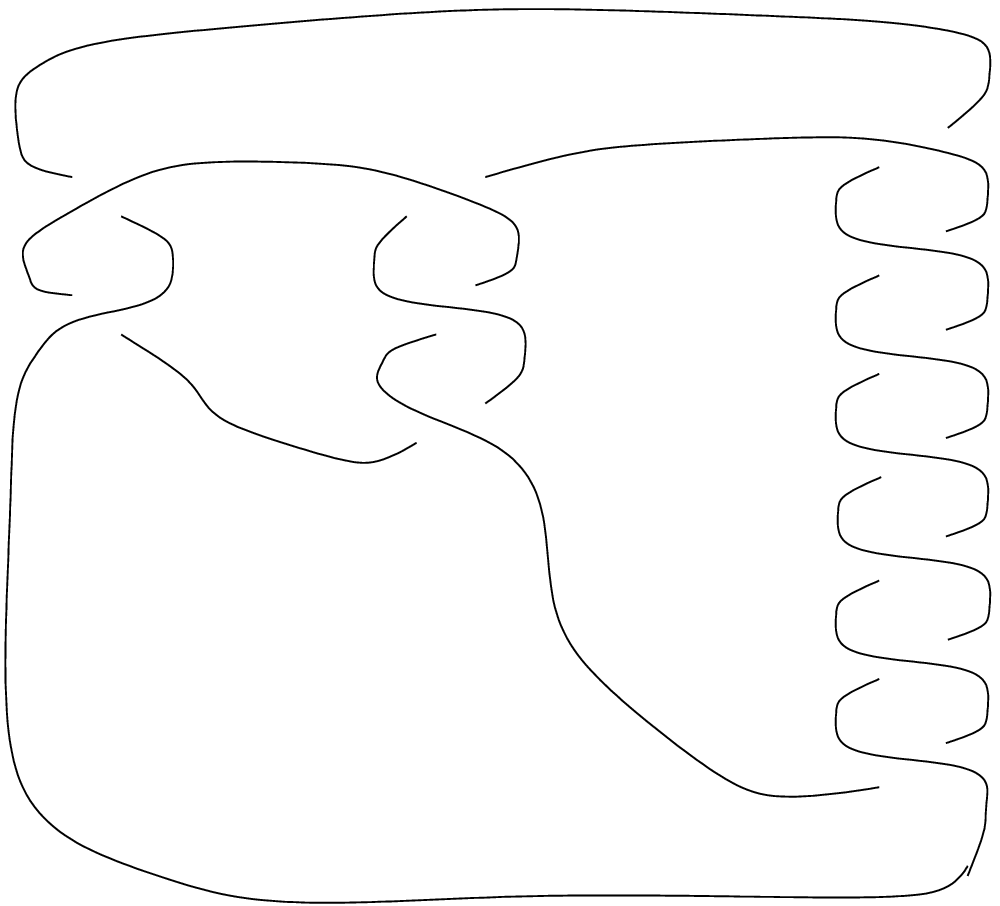,height=1in}}

The $(-2,3,7)$ pretzel knot is equivalent to the $(2,3,7,-1)$ pretzel
knot.  Let $K_{p_1,\dots,p_k}$ be the $(p_1,\dots,p_k,-1,\dots,-1)$-pretzel
link, where the number of ``-1"s is $k-2$.  
It turns out that the Alexander polynomial of $K_{p_1,\dots,p_k}$ is related
to the denominator of the growth series of $T_{p_1,\dots,p_k}$ 
\cite{Hironaka98}.

\begin{theorem} The pretzel link $K_{p_1,\dots,p_k}$
is fibered and has Alexander polynomial 
$$
\Delta_{p_1,\dots,p_k}(-x).
$$
\end{theorem}

Thus, Theorem~\ref{delta-theorem} implies the following.

\begin{corollary}\label{pretzel-corollary} Among
pretzel links $K_{p_1,\dots,p_k}$, the Mahler measure of the
Alexander polynomial is minimized by $K_{2,3,7}$.
\end{corollary}

Although the relation between the Alexander polynomial
and the denominator of the growth series has not been fully explained
there is a natural relation between the pretzel links and the polygonal
reflection groups which we describe in Section~\ref{pairing-section}.

\subsection{$E_{10}$ diagram}

The Coxeter $E_{10}$ diagram can be thought 
of as the $(2,3,7)$-star Coxeter graph as can be seen by comparing
Figure~\ref{E10} with Figure~\ref{star}. 
\makefig{$E_{10}$-Coxeter graph}{E10}
{\psfig{figure=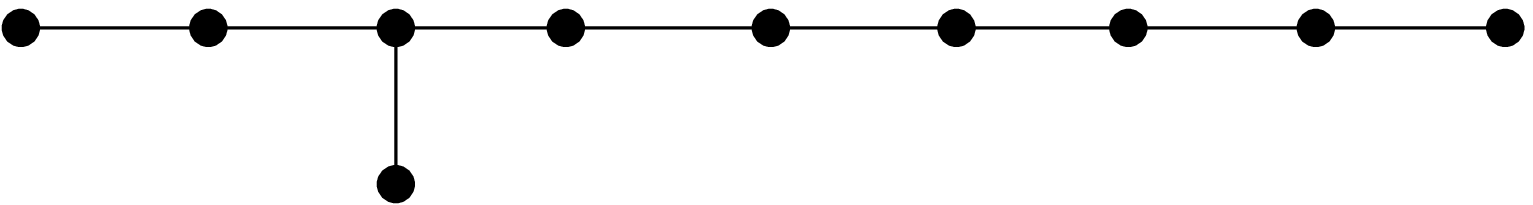,height=0.3in}}

McMullen observes \cite{McMullen:Coxeter} that the
characteristic polynomial of the Coxeter element is Lehmer's
polynomial $p_L(x)$, and its leading eigenvalue is $\alpha_L$.
Furthermore, he shows the following.

\begin{theorem}  Let $C$ be the Coxeter element of a Coxeter system,
and let $\lambda(C)$ be the spectral radius of $C$.  Then 
among non-spherical and non-affine Coxeter systems, $\lambda(C)$
achieves its minimum when $(G,S)$ is the Coxeter system corresponding
to the $E_{10}$ diagram.  
\end{theorem}

This solves Lehmer's problem for Coxeter elements of Coxeter systems,
and as a consequence also for the monodromy of Coxeter links, generalizing
Corollary~\ref{pretzel-corollary}.

\section{Correspondence between stars and polygons}\label{pairing-section}

The results in the previous sections show that the growth 
rate of $(G,S)$ corresponding to the Coxeter system
$T_{p_1,\dots,p_k}$ and the spectral radius of the Coxeter
system $\Star(p_1,\dots,p_k)$ are equal.

\makefig{$\Star(p_1,\dots,p_k)$ and $T_{p_1,\dots,p_k}$}{dual}
{\psfig{figure=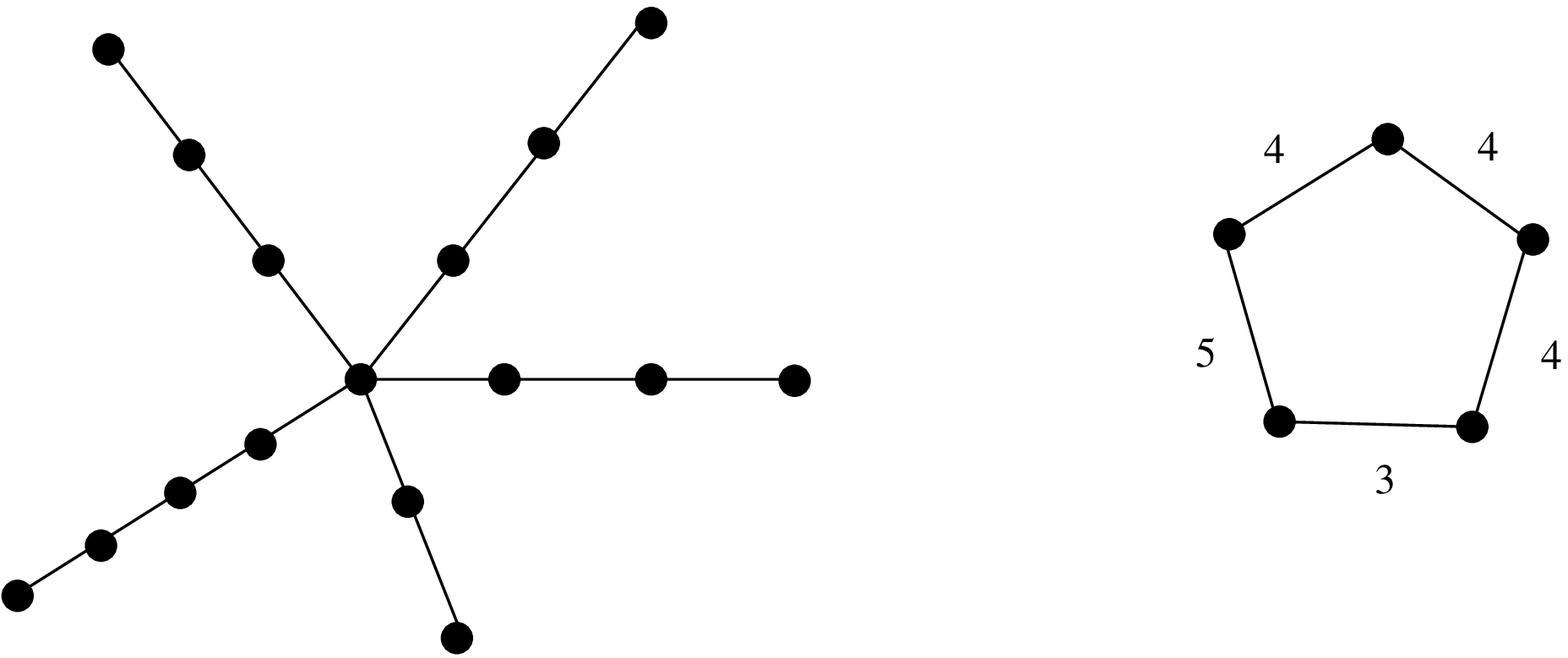,height=1.25in}}

Here we give a topological relation between the polygonal- and
star-Coxeter systems.
Let $P$ be the sphere $S^2$, the Euclidean plane $\E^2$, or the
hyperbolic plane $\H^2$.  Let $T(P)$ be the unit tangent bundle
of $P$.  Let $G$ be the $p_1,\dots,p_k$-polygonal reflection group 
acting on $P$.

It is not hard to see that the
action of $G$ on $T(P)$ makes $T(P)$ a branched cover over 
$S^3$ with branching index 2 on a link $K$.   

Let $X$ be a manifold and
$Y$ a formal sum of codimension one disjoint submanifolds 
of $X$ with integer coefficients. 
We denote by $(X;Y)$ the orbifold $X$ with orbifold 
singularities at $Y$ of the appropriate multiplicity.

The covering $T(P) \rightarrow S^3$ factors through
an unbranched covering of a manifold $M$ which double covers $S^3$
as in the following diagram:    
$$
\begin{CD}
T(P) \\
@V{G^{(2)}}VV \\
M @>>> (S^2;p_1[x_1] + \cdots + p_k[x_k])\\
@V{\Z_2}VV \\
(S^3;2[K])
\end{CD}
$$

A double covering of $S^3$ branched 
along any Montesinos link $K$ is
Seifert fibered over $S^2$ (see, for example, \cite{B-Z:Knots}.) 
The $p_1,\dots,p_k$-star Coxeter link 
is only a particular example of a Montesinos link
for which the double cover is Seifert fibered over
$(S^2;p_1[x_1] + \cdots + p_k[x_k])$, but there are many others
with different Alexander polynomials.
For the particular 3-manifold $M$ defined above, however, 
one can show the following is true.

\begin{proposition}\label{pairing-proposition} The link $K$ 
coming from the quotient of the homogeneous space $T(P)$ by 
the $p_1,\dots,p_k$-polygonal Coxeter group is the Coxeter 
link associated to the $p_1,\dots,p_k$-star Coxeter system. 
\end{proposition}

\subsection{Example: Klein Singularities}

The above discussion generalizes a phenomenon encountered in the
theory of isolated hypersurface singularities.

Let $\widetilde{G}$ be a finite subgroup of 
$\mbox{SU}(2)$ acting in the usual way on $\C^2$.  The 
quotient $X$ is a hypersurface in $\C^3$ with isolated singularity
called a Klein singularity.  
The resolution diagrams of these singularities are exactly the
Coxeter graphs (or Dynkin diagrams) of the simply laced spherical
Coxeter systems.  The correspondence between
finite subgroups of $\mbox{SU}(2)$ and the Dynkin diagrams is
known as McKay's Correspondence \cite{McKay:Correspondence} (see
also \cite{G-V:McKay} and \cite{Slodowy:Klein}.)  

Coxeter links can be used to give a topological description of 
the correspondence.
Let $\widetilde{G}$ be a finite subgroup of $\mbox{SU}(2)$.
Then $\widetilde{G}$ is the binary
extension of a finite subgroup $G^{(2)}$ of $\mbox{SO}(3)$ 
which in turn (with the exceptional case of $G^{(2)} = \Z_n$) 
is an index 2 subgroup of a $(p,q,r)$-triangle group $G$.  
The group $G^{(2)} = \Z_n$ corresponds to the $(2,n)$-triangle
group.
Relating each group with the $(p,q,r)$-star diagram gives the
McKay correspondence, as seen in the following diagram.

\begin{center}
\begin{tabular}{|  c | c | c |}
\hline
${\Z}_n$ & $T_{2,n}$ &$A_{n+1}$\\
\hline
${\cal D}_n$ & $T_{2,2,n}$  &$D_{n+2}$\\
\hline
$S_3$ & $T_{2,3,3}$  &$E_6$\\
\hline
$S_4$ & $T_{2,3,4}$  &$E_7$\\
\hline
$A_5$ & $T_{2,3,5}$ &$E_8$\\
\hline
\end{tabular}
\end{center}

For these triples $p,q,r$ in the above diagram,
the $(p,q,r)$-star Coxeter link $K$ is algebraic,
and the induced branched covering of the link $M$ of 
the singularity $X$ over $S^3$ is induced by a
generic projection of $X$ to $\C^2$.  Thus, we obtain
the following commutative diagram. 
$$
\begin{CD}
\R\P^3 = T(S^2) @<{2}<< S^3 \subset \C^2\\
@V{G^{(2)}}VV @V{\widetilde G}VV\\
M @>\subset>> X\\
@V{\Z_2}VV @V{\rho}VV\\
S^3 @>\subset>> \C^2
\end{CD}
$$
A minimal desingularization of $X$ can be obtained by
blowing up $\C^2$ to desingularize the branch curve
of the map $\rho$, and pulling back over $X$. 
$$
\begin{CD}
X @<{\widetilde \sigma}<< \widetilde X\\ 
@V{\rho}VV @V{\widetilde {\rho}}VV\\
\C^2  @<{\sigma}<< \widetilde {\C^2}
\end{CD}
$$
One way to verify that the resolution diagram for $X$ 
must equal the Dynkin diagram corresponding to the link,
is to note that the effect of desingularizing the branch curve 
is the same as unknotting the branch link by a sequence
of simple Dehn twists on suitable closed curves in the 
complement. 
\newpage

\bibliographystyle{math}
\bibliography{math}

\noindent
Eriko Hironaka\\
Department of Mathematics\\
Florida State University\\
Tallahassee, FL 32306\\
Email: hironaka@math.fsu.edu

\end{document}